\documentclass[11pt]{amsart}
\usepackage{cite}
\usepackage{pinlabel}
\usepackage{amsmath} 
\usepackage{graphicx} 
\setkeys{Gin}{keepaspectratio}
\usepackage{hyperref}
\usepackage{url}
\usepackage{bm}
\usepackage{mathabx}
\usepackage[left=1.25in,top=1in,right=1.25in,bottom=1in,head=.1in]{geometry}
\usepackage{xcolor}
\usepackage{tikz}
\usepackage{adjustbox}
\usepackage[normalem]{ulem} % 
\usepackage{enumitem}

\usepackage{verbatim}
\newcommand{\items}{\begin{itemize}[leftmargin=25pt,rightmargin=15pt]
  \setlength\itemsep{2pt}}
\newcommand{\stopitems}{\end{itemize}}

\usepackage[textsize=scriptsize]{todonotes}

\usepackage{fancyhdr}
\newtheorem{theorem}{Theorem}[section] % numbered theorems, lemmas, etc
\newtheorem*{theorem*}{Theorem}
\newtheorem{lemma}[theorem]{Lemma}

\newtheorem*{conjecture*}{Conjecture}
\newtheorem*{question*}{Question}
\newtheorem*{lemma*}{Lemma}
\newtheorem{corollary}[theorem]{Corollary}
\newtheorem*{corollary*}{Corollary}
\theoremstyle{definition}

\newtheorem{remark}[theorem]{Remark}

\newtheorem*{example*}{Example}
\newtheorem*{remark*}{Remark}
\newtheorem*{remarks*}{Remarks}
\newtheorem*{addenda*}{Addenda}
\newtheorem*{construction*}{Construction}

\newcommand{\sss}{S^2 \times S^2} % \ss is German double-s

\newcommand{\Z}{\mathbb Z}
\newcommand{\bc}{\mathbb C}

\newcommand{\calo}{\mathcal O}

\newcommand{\id}{\textup{id}}

\renewcommand{\phi}{\varphi}

\DeclareMathOperator{\diffp}{\it Diff^+}

\newcommand{\cs}{\mathbin{\#}}
\newcommand{\cpone}{\bc P^1}

\newcommand{\cptwo}{\bc P^2}

\newcommand{\cptwobar}{\smash{\overline{\bc P}^2}}

\newcommand{\unred}[1]{ \ignorespaces}  %% to omit that red text
\title{Wall's stable realization for diffeomorphisms of definite 4-manifolds\\ \today}
\author[Daniel Ruberman]{Daniel Ruberman}
\address{Department of Mathematics, MS 050\newline\indent Brandeis
University \newline\indent Waltham, MA 02454}
\email{ruberman@brandeis.edu}
\author[Sa\v{s}o Strle]{Sa\v{s}o Strle}
\address{Faculty of Mathematics and Physics\newline\indent University of Ljubljana\newline\indent Jadranska 19\newline\indent 1000 Ljubljana\newline\indent Slovenia}
\email{saso.strle@fmf.uni-lj.si}
\thanks{We thank Anubhav Mukherjee for his observations about Wall's work that led to this note. The first author was partially supported by NSF Grant DMS-1928930 while he was in residence at the Simons Laufer Mathematical Sciences Institute (formerly known as MSRI), as well as NSF FRG Grant DMS-1952790. The second author was partially supported by Slovenian Research Agency (ARRS) Research program P1-0288.\\
Math.~Subj.~Class.~2020: 57K40 (primary), 57R50 (secondary).}

\begin{document}
\setlength{\headheight}{12.0pt}
\begin{abstract}
Let $X$ be a smooth simply connected closed $4$-manifold with definite intersection form. We show that any automorphism of the intersection form of $X$ is realized by a diffeomorphism of $X \cs \sss$. This extends and completes Wall's foundational result from 1964.
\end{abstract}
\maketitle
\section{Introduction}
In 1964, C.T.C.~Wall published two foundational papers~\cite{wall:diffeomorphisms,wall:4-manifolds} about closed simply connected $4$-manifolds. The main theorem of~\cite{wall:diffeomorphisms} is a `1-stable realization result' stating that for any smooth closed simply connected manifold (henceforth denoted $X$), every automorphism of the intersection form of $X \cs \sss$  is realized by a diffeomorphism, provided that the intersection form of $X$ is either indefinite or of rank at most $8$. Wall used this to show in~\cite{wall:4-manifolds} that homotopy equivalent manifolds are h-cobordant, and therefore diffeomorphic after connected sum with some number of copies of $\sss$.  Although he did not explicitly consider the topological case, his results in the smooth case were crucial in Freedman's topological classification~\cite{disc,freedman,freedman-quinn} some 20 years later. One consequence of Freedman's work in the topological case is that one does not have to add any $\sss$ summands: any automorphism of the intersection form of $X$ is realized by a self-homeomorphism of $X$. 

The restriction on the intersection form is necessary, as pointed out by Anubhav Mukherjee. Friedman and Morgan~\cite{friedman-morgan:complex-surfaces-I,friedman-morgan:complex-surfaces-II,friedman-morgan:sw-surfaces} showed that not every automorphism of the intersection form of $\cptwo  \cs^k \cptwobar$ for $k \geq 10$ is realized by a diffeomorphism.  Since $\cptwo \cs^k \cptwobar \cong  \cs^{k-1} \cptwobar \cs \sss$, this shows that Wall's 1-stable realization result fails for definite forms of rank $9$ and higher.

In this note we show that a portion of the 1-stable realization result holds for \emph{all} smooth simply connected closed $4$-manifolds, without restriction on the rank of the intersection form.  Our result (stated more precisely below) says that any automorphism of the intersection form of $X$ is realized on $X \cs \sss$. In particular, the automorphisms considered by Friedman and Morgan are not of this form. 

\section{Results}
Let $X$ be a smooth simply connected 4-manifold.  Choosing an orientation on $X$, we denote by $Q_X$ the intersection form on $H_2(X)$, and by $\calo(Q_X)$ its automorphism group. Passing from an orientation-preserving self-diffeomorphism of $X$ to its induced map on $H_2(X)$ gives a natural homomorphism from $\diffp(X) \to \calo(Q_X)$, and we say that $A \in \calo(Q_X)$ is realized by a diffeomorphism if $A = f_*$ for some $f\in \diffp(X)$. There is also a natural inclusion $A \to A \oplus \id_{H}$ of $\calo(Q_X)$ into $\calo(Q_{X\cs \sss})$, where $H=Q_{\sss}$ is the hyperbolic form. It is straightforward to see that for any $f \in \diffp(X)$ there is a diffeomorphism of $X\cs \sss$ restricting to the identity on $\sss$ that realizes $f_* \oplus \id_{H}$.
\begin{theorem}\label{T:realize}
Let $X$ be a smooth simply connected closed $4$-manifold and let $A \in \calo(Q_X)$.  Then $A \oplus \id_{H}$ is realized by a diffeomorphism of $X \cs \sss$.
\end{theorem}

\begin{proof}
Since Wall has shown the result for indefinite manifolds, it suffices to work with a definite manifold, which we take to be negative definite. By Donaldson's diagonalization theorem~\cite{donaldson}, we can identify $Q_X$ with $L_k = \langle -1 \rangle^k$, with a basis $\{e_1,\ldots, e_k\}$ for which $Q_X(e_i, e_j) = -\delta_{ij}$. We follow standard usage and write $\alpha \cdot \beta$ for $Q_X(\alpha,\beta)$ and refer to $Q_X(\alpha,\alpha)=\alpha \cdot \alpha$ as the \emph{square} $\alpha^2$ of $\alpha$. We now make use of some results of Wall's papers~\cite{wall:automorphisms-II}  and~\cite{wall:diffeomorphisms}.  The first is the observation from~\cite{wall:automorphisms-II} that $\calo(L_k)$ is generated by permutations of the basis vectors $e_i$, together with reflections $S(e_i)$ defined by 
\[
S(e_i)(e_j) = 
\begin{cases}
-e_j \quad \text{if }j=i, \text{and}\\
e_j  \quad \text{otherwise.}
\end{cases}
\]
The second is the construction of certain diffeomorphisms $f_w^i$ of $X \cs \sss$ realizing the automorphisms $E_w^i \in \calo(Q_{X\cs \sss}) \cong \calo(Q_X \oplus H)$ defined by Eichler~\cite{eichler:quadratic}.  Here $w$ is an element of $H_2(X)$ with \emph{even} square $w^2=2s$ and $i = x, y$, where $x$ and $y$ are the standard generators of $H_2(\sss)$ satisfying $x^2=y^2=0$ and $x\cdot y=1$. Then
\[
E_w^y(u)=u-(u \cdot w)y,\quad E_w^y(x)=x+w-sy,\quad E_w^y(y)=y
\]
for $u \in H_2(X)$. $E_w^x$ is defined similarly with the roles of $x$ and $y$ reversed.

In Wall's proof, for an indefinite form or for a definite form of rank at most 8 that is of even type, one sees that the automorphism group $\calo(Q_X)$ is contained in the subgroup generated by the diffeomorphisms $f_w^i$ and diffeomorphisms of $\sss$.
For $Q_X$ odd, which is the situation in our case, one additional diffeomorphism $\phi$ is required to arrive at the same conclusion.  This can be described by choosing (as in~\cite[Theorem 1]{wall:diffeomorphisms}) a diffeomorphism
\[
g_w: X \cs \sss \overset{\cong}{\longrightarrow} X\cs \cptwo \cs \cptwobar
\] 
that depends on a vector $w \in H_2(X)$ of \emph{odd} square  
$w^2=r=2s-1$.  This is analogous to $f_w^i$ above; its action on homology is given by
\[
G_w(u)=u-(u \cdot w)(h-k),\quad G_w(x)=w+h-s(h-k),\quad G_w(y)=h-k,
\]
where $h$ (resp. $k$) is the generator of $H_2(\cptwo)$ (resp. $H_2(\cptwobar)$) corresponding to $\cpone$.
Then $\phi_w$ is obtained by using $g_w$ to transport the diffeomorphism given by complex conjugation on the $\cptwo$ summand back to $X\cs \sss$.  The induced automorphism $\Phi_w$ is then given by
\begin{align*}
\Phi_w(u) & =u+(u \cdot w)(2x+(1-r)y-2w)\\ 
\Phi_w(x) & =(1-r)w+rx-(r-1)^2y/2\\ 
\Phi_w(y) & =2w-2x+ry
\end{align*}
where $u, x$ and $y$ are as above.

Suppose now that $Q_X$ is $L_k$ for some $k$. Then for any element $e\in H_2(X)$ of square $-1$ 
\[
S(e)=E_{2e}^x \circ C_* \circ E_{2e}^y \circ \Phi_e,
\]
where $C: \sss \to \sss$ is complex conjugation in both factors, may be realized by a diffeomorphism.  Similarly for a pair $e_1, e_2$ of linearly independent elements of square $-1$ their transposition may be expressed as
\[
T_* \circ  E_{e_1-e_2}^y \circ E_{e_2-e_1}^x \circ E_{e_1-e_2}^y,
\]
where $T: \sss \to \sss$ interchanges the factors, showing it is also realized.  Thus all of the generators of  $\calo(Q_{X})$ are realized by diffeomorphisms of $X \cs \sss$ that preserve the $Q_X$ summand.
\end{proof}

As Wall remarks in~\cite{wall:diffeomorphisms}, his 1-stable realization theorem holds for $4$-manifolds with homology sphere boundary (with the same restriction on intersection forms as in the closed case).  
% Let $\calw$ be the collection of all unimodular forms for which Wall's results hold.
Let $E_8$ be the unique negative definite even unimodular lattice of rank 8.
As in Theorem \ref{T:realize} we obtain the following slightly more general version of 1-realization in this case. 
\begin{corollary}\label{C:realize-boundary-sphere}
Let $X$ be a simply connected $4$-manifold with boundary an integral homology sphere $Y$. Suppose that the intersection form $Q_X$ is either indefinite, diagonalizable or the direct sum of $\pm E_8$ and a diagonalizable form. Then for any $A \in \calo(Q_X)$  there exists a diffeomorphism of $X \cs \sss$ realizing $A \oplus \id_{H}$ and fixing the boundary $Y$.
\end{corollary}
\begin{proof}
For $Q_X$ indefinite or of rank at most 8 this is Wall's theorem. For diagonalizable $Q_X$ the same construction as in the proof of the previous theorem applies. Suppose now $Q_X=E_8 \oplus L_k$.  Then any automorphism of $Q_X$ respects the splitting so again can be realized.

Note that all diffeomorphisms involved in realizations may be chosen so that they fix $Y$.  Recall from~\cite{wall:diffeomorphisms} that $f_w^i$  is on $X$ supported in a neighborhood of a representative of the class $w$ which may be chosen disjoint from $Y$. The same applies to $\phi_w$. 
\end{proof}

Note that the restriction of the corollary on definite $Q_X$ holds if $Y$ admits an orientation reversing diffeomorphism.
\begin{lemma}\label{L:reverse-orientation}
Let $X$ be a simply connected $4$-manifold with boundary an integral homology sphere $Y$. If the intersection form $Q_X$ is definite and $Y$ admits an orientation reversing diffeomorphism, then $Q_X$ is diagonalizable. 
\end{lemma}
\begin{proof}
Let $g:Y \to Y$ be an orientation reversing diffeomorphism. Then $Z=X \cup_g X$ is a closed definite manifold with $Q_Z=Q_X \oplus Q_X$. Since by Donaldson's diagonalization theorem $Q_Z$ is diagonalizable, so is $Q_X$. Indeed, supposing $X$ is positive definite there is a basis $\{e_1,\ldots, e_{2k}\}$ of $H_2(Z)=H_2(X) \oplus H_2(X)$ for which $Q_Z(e_i, e_j) = \delta_{ij}$. Since a vector of square $1$ cannot be the sum of two orthogonal vectors, each $e_i$ belongs to one copy of $H_2(X)$. 
\end{proof}

It is an important question to what extent above results generalize to 4-manifolds with other boundaries.  Let $X$ be a simply connected 4-manifold with boundary a rational homology sphere $Y$.  Denote by $(L,Q)=(H_2(X),Q_X)$ the intersection lattice of $X$ and by $V$ the real vector space underlying $L$; note that $Q$ extends uniquely to $V$. Since the 2-dimensional homology groups in the exact sequence
\[
0 \to H_2(X) \to H_2(X,Y) \to H_1(Y) \to 0
\]
are free abelian and $H_1(Y)$ is finite, it follows that the dual lattice $L^*=H_2(X,Y)$ is a sublattice of $V$. Let $A$ be an automorphism of $(L,Q)$, so 
\[
Q(Ax,Ay)=Q(x,y)
\]
for all $x, y \in L$ and hence also for all $x, y \in V$. Then for any $w \in L^*$, it follows that $Aw \in L^*$ since for all $x \in L$
\[
Q(Aw,x)=Q(Aw,Az)=Q(w,z) \in \Z,
\]
where $z \in L$ is determined by $x=Az$.  Hence $A$ extends to an automorphism of $L^*$. In particular, any automorphism $A \in \calo(Q_X)$ extends to an automorphism of $H_2(X,Y)$ and through this defines an action on $H_1(Y)$.  We say that $A$ \emph{acts trivially} on $Y$ if this action is trivial.  
\begin{theorem}\label{T:realize-boundary}
Let $X$ be a simply connected $4$-manifold with boundary a homology lens space $Y$ of order $d>1$, i.e. $H_1(Y) \cong \Z/d\Z$.  Suppose $Q_X \cong Q \oplus \langle d \rangle$ with $Q$ as in Corollary \ref{C:realize-boundary-sphere}.  Then for any $A \in \calo(Q_X)$ acting trivially on $Y$ there exists a diffeomorphism of $X \cs \sss$ realizing $A \oplus \id_{H}$ and fixing the boundary $Y$.
\end{theorem}
\begin{proof}
Let $z \in H_2(X)$ be a generator supporting $\langle d \rangle$. Then $H_2=U \oplus \Z z$ with $U$ supporting $Q$. Since $z$ is orthogonal to $U$ and of square greater than 1, $A$ respects the above splitting. In particular, $A(z)$ is $\pm z$. The image of $z$ in $H_2(X,Y)$ is of the form $dw$ and the image of $w$ generates $H_1(Y)$.  For $d>2$ the condition that the action of $A$ on $Y$ is trivial implies $A(z)=z$, whereas for $d=2$ both the possibilities may occur. 

If $A(z)=z$ then $A$ is realized by the same argument as before. If $A(z)=-z$ then the action of $A$ on $z$ can be in $\calo(Q_X \oplus H)$ expressed as
\[
(CT)_* \circ E_z^y \circ E_z^x  \circ E_z^y
\]
which is realized by a diffeomorphism fixing the boundary.
\end{proof}

\section{Discussion}
We close with a few comments regarding realization problems.
\begin{remark} The diffeomorphisms $f^i_w$ and $\phi_w$ in Wall's paper are also described in terms of handlebody pictures in Kirby's  book~\cite[Chapter X]{kirby:4-manifolds}. Kirby assumes that $X$ has a handlebody decomposition with no $1$-handles, but the arguments go through virtually unchanged in the presence of $1$-handles. The main observation is that (perhaps after changing the handle structure by standard operations) a primitive homology class $w$ will be represented by a single $2$-handle whose attaching map goes over each $1$-handle algebraically $0$ times.
\end{remark}
\begin{remark}
For the simplest negative definite manifold $\cs^k \cptwobar$, the generators are realized (without any stabilizations) by permuting the $\cptwobar$ summands or by complex conjugation on the $i^{th}$ summand.  There are no known examples of exotic definite manifolds, so it is possible in principle that our theorem could be deduced by classifying definite manifolds up to diffeomorphism. But this seems rather difficult at present. 
\end{remark}
\begin{remark} It is an interesting question as to whether the result of Corollary \ref{C:realize-boundary-sphere} holds for arbitrary intersection forms, which of course can be $Q_X$ for a manifold $X$ with non-empty homology sphere boundary.  To use the current argument, one would have to solve the algebraic problem of whether for an arbitrary form $Q$,  the subgroup of $\calo(Q \oplus H)$ generated by  the $E^i_w$, elements of $\calo(H)$, and the automorphisms $\Phi_w$ actually contains all automorphisms of $Q$.  
\end{remark}
\begin{remark}
Wall's papers show, broadly speaking, that classification theorems in smooth $4$-manifold theory are possible if one works with stabilized objects, where stabilization refers to connected sum with $\sss$. Theorem~\ref{T:realize} completes Wall's result that for realization results, one stabilization suffices. There are no known examples of exotic pairs of closed manifolds that require more than one stabilization to become diffeomorphic; see for instance~\cite{mandelbaum:decomp,mandelbaum:irrational,mandelbaum-moishezon:algebraic,akbulut:fs-knot-surgery,auckly:stable,baykur:stable}, making `one stabilization suffices' an attractive conjecture.  However, recent work of Kang \cite{kang} shows that one stabilization is not enough for manifolds with boundary.

For other types of smooth phenomena, the situation is more complicated. The papers~\cite{baykur-sunukjian:round,auckly-kim-melvin-ruberman:isotopy,auckly-kim-melvin-ruberman-schwartz:1-stable} show that many families of closed mutually exotic surfaces are isotopic after one stabilization; this yields exotic diffeomorphisms~\cite{ruberman:isotopy} and higher dimensional families of diffeomorphisms~\cite{auckly-ruberman:families} that become standard after one stabilization. However, an example of Jianfeng Lin~\cite{lin:dehn} shows that one stabilization is not always sufficient to provide an isotopy, and examples of Lin and Mukherjee~\cite{lin-mukherjee:surfaces} show that a single stabilization is not sufficient when working with surfaces with boundary and diffeomorphisms on manifolds with boundary. Finally (for the present) forthcoming work of Konno-Mukherjee-Taniguchi~\cite{konno-mukherjee-taniguchi:codim1} shows that there are exotically embedded $3$-manifolds which \emph{no} amount of stabilization can make standard.
\end{remark}

%%%%%%%%% bib

\def\cprime{$'$}
\providecommand{\bysame}{\leavevmode\hbox to3em{\hrulefill}\thinspace}

\end{document}